# ON SOME PROBLEMS IN THE ARTICLE EFFICIENT LIKELIHOOD ESTIMATION IN STATE SPACE MODELS BY CHENG-DER FUH
## [ANN. STATIST. 34 (2006) 2026–2068]

By Jens Ledet Jensen

*University of Aarhus*

**1. Introduction.** Upon reading the paper *Efficient Likelihood Estimation in State Space Models* by Cheng-Der Fuh I found a number of problems in the formulations and a number of mathematical errors. Together, these findings cast doubt on the validity of the main results in their present formulation. A reformulation and new proofs seem quite involved.

The paper, *Efficient Likelihood Estimation in State Space Models* deals with asymptotic properties of the maximum likelihood estimate in hidden Markov models. The hidden Markov chain is $X_n$, and the observed process is $\xi_n$ where $\xi_n$ conditioned on the past and the hidden process depends on $(X_n, \xi_{n-1})$ only. The approach used is to add an iterated function system $M_n$, and to consider the Markov process $(X_n, \xi_n, M_n)$. This is very much akin to the method in Douc and Matias [1], and I will use this article as a background for my comments.

**2. Problems.**

2.1. *Definition of iterated function system.* The first basic definition in the paper is a function $\mathbf{P}_\theta(\xi_j): \mathbf{M} \to \mathbf{M}$ that maps a function $h \in \mathbf{M}$ into a new function in $\mathbf{M}$ (page 2031),

$$\mathbf{P}_\theta(\xi_j)h(x) = \int_{y \in \mathcal{X}} p_\theta(x,y)f(\xi_j;\theta|y,\xi_{j-1})h(y)m(dy).$$

[It is unclear why the author states that $\mathbf{P}_\theta(\xi_j)$ is a function on $(\mathcal{X} \times \mathbf{R}^d) \times \mathbf{M}$ where $\mathcal{X}$ is the state space of the Markov chain.] The paper next defines the composition $\mathbf{P}_\theta(\xi_{j+1}) \circ \mathbf{P}_\theta(\xi_j)h$ by first applying $\mathbf{P}_\theta(\xi_{j+1})$ to $h$ and then applying $\mathbf{P}_\theta(\xi_j)$ to the result. Using these two definitions we have

$$\mathbf{P}_\theta(\xi_n) \circ \cdots \circ \mathbf{P}_\theta(\xi_1) \circ \mathbf{P}_\theta(\xi_0)\pi_\theta$$
$$= \int \pi_\theta(x_n) \left\{ \prod_{j=n}^{1} p_\theta(x_{j-1}, x_j) f(\xi_j; \theta | x_j, \xi_{j-1}) m(dx_j) \right\} f(\xi_0; \theta | x_0) m(dx_0).$$

The argument presented in the paper then appears to assume that this expression depends on some $x$ and performs an integration before claiming that the result is the joint density $p_n(\xi_0, \ldots, \xi_n; \theta)$. This is clearly not correct since $\pi_\theta(x_n)$ appears in the expression instead of $\pi_\theta(x_0)$.

Following the work of Douc and Matias [1] one would instead use the definition

(1) $$\mathbf{P}_\theta(\xi_j)h(x) = \int_{y \in \mathcal{X}} p_\theta(y, x) f(\xi_j; \theta | y, \xi_{j-1}) h(y) m(dy);$$

that is, the integration is with respect to the first variable in $p_\theta(y, x)$ instead of the second. Changing the definition of $\mathbf{P}_\theta(\xi_0)$ correspondingly and using ordinary composition of functions, one finds that $p_n(\xi_0, \ldots, \xi_n; \theta)$ equals the integral of $\mathbf{P}_\theta(\xi_n) \circ \cdots \circ \mathbf{P}_\theta(\xi_1) \circ \mathbf{P}_\theta(\xi_0)\pi_\theta$ with respect to $x_{n+1}$. However, making this change necessitates a new proof for the first part of Lemma 3 on page 2056. Comparing with Douc and Matias ([1], Proposition 1) we see that this is one of the places where the latter authors use the stronger assumptions of that paper on the Markov chain.

Turning to the iterated function system, Fuh's paper defines this as

$$M_n = \mathbf{P}_\theta(\xi_n) \circ \cdots \circ \mathbf{P}_\theta(\xi_1) \circ \mathbf{P}_\theta(\xi_0)$$

[formula (5.6), page 2045]. Taking this literally, and using the definitions in Fuh's paper, this is actually a mapping that takes a function as input and turns it into a constant. Instead $M_n$ should be a function obtained by applying a mapping to $M_{n-1}$. This is achieved when using the definition suggested in (1) and adding $\pi_\theta$ to the right-hand side of $M_n$ above.

2.2. *Harris recurrence of iterated function.* Whether or not we make the changes suggested in the previous subsection, $M_n$, defined on page 2045, is related to the density of $(\xi_0, \ldots, \xi_n)$. Making the change suggested in (1) above we have precisely $M_n(x_{n+1}) = p(x_{n+1}, \xi_0, \ldots, \xi_n)$. Such an expression will typically tend to either zero or infinity. However, in Lemma 4 on page 2046 Fuh claims that $(X_n, \xi_n, M_n)$ is a Harris recurrent Markov chain. It is difficult to pinpoint the exact origin of this problem. The Harris recurrence is established in Lemma 4 which in its formulation uses a measure $Q$ from Theorem 1 (in the formulation there are two $Q$'s, but these are different).

So we need to establish Theorem 1 before proving Lemma 4. In Lemma 3 it is stated that the Markov iterated function system satisfies Assumption K. In Remark 1 (page 2035) Fuh says that Assumption K is different from the assumptions of Theorem 1. He then goes on to say that if Assumption K is supplemented with the extra assumption that $(Y_n, M_n)$ is a Harris recurrent Markov chain, then Theorem 1 still holds. This, therefore, seemingly looks like a circular argument.

Comparing again with Douc and Matias [1] they consider instead $M_n(x_{n+1}) = p(x_{n+1}|\xi_0, \ldots, \xi_n)$. However, if we make this change we have introduced a new iterated function system, and a revised version of Lemma 3 is needed which presumably will lead to a different set of assumptions.

2.3. *Asymptotic properties of score function and observed information.* The asymptotic normality of the score function is stated in Lemma 6 (page 2048). In the proof of Lemma 6 (page 2060) the author appeals to Corollary 1. The latter gives a central limit theorem for a sum of the form $\sum_{j=1}^n g(M_j)$. However, the paper wants to use this result on the sum $\sum_{j=1}^n \frac{\partial}{\partial \theta} g(M_{j-1}, M_j)$. This looks innocent, but since $\theta$ appears in the iteration of $M_n$ this is not on the form $\sum_{j=1}^n \tilde{g}(M_{j-1}, M_j)$. Instead one needs to consider a new iterated function system. This is what is done in Appendix D of Douc and Matias [1].

Similarly, it is stated that the proof of the main Theorem 5 follows a standard argument. However, comparing with Douc and Matias [1] (Appendix D.3) it seems that yet another iterated function system is needed to deal with the convergence of the observed information.

2.4. *Generality of conditions.* Assumption C5 on page 2043 restricts the dependency of the observed process on the hidden process. For the example considered in (b) on page 2044 one needs to consider

$$\sup_{y,z\in\mathcal{X}} \frac{f(\xi_0;\theta|y)f(\xi_1;\theta|y,\xi_0)}{f(\xi_0;\theta|z)f(\xi_1;\theta|z,\xi_0)}$$

$$= \sup_{y,z\in\mathcal{X}} \frac{\exp\{-1/2(\xi_0-y)^2 - 1/2(\xi_1-y)^2\}}{\exp\{-1/2(\xi_0-z)^2 - 1/2(\xi_1-z)^2\}}$$

$$= \sup_{y,z\in\mathcal{X}} \exp\{z^2 - y^2 + (\xi_0+\xi_1)(y-z)\} = \infty.$$

Thus C5 is not satisfied (this seems to be contrary to the claim on page 2054 line 8 from the bottom).

## REFERENCES

[1] DOUC, R. and MATIAS, C. (2001). Asymptotics of the maximum likelihood estimator for general hidden Markov models. *Bernoulli* **7** 381–420. MR1836737

Department of Mathematical Sciences
University of Aarhus
DK-8000 Aarhus C
Denmark
E-mail: jlj@imf.au.dk